\documentclass[11pt,leqno,a4paper,twoside,final]{article}
\usepackage[T1]{fontenc}
\usepackage{exscale,ifthen,array,amsmath,amscd,amssymb,amsthm,natbib}

\allowdisplaybreaks[1]

\DeclareMathOperator{\Id}{Id}

\DeclareMathOperator{\intr}{int}
\DeclareMathOperator{\Prb}{\mathrm{P}}
\DeclareMathOperator{\Mean}{\mathbf{E}}

\bibpunct[:\,]{(}{)}{,}{a}{,}{,}
\begin{document}


\newcommand{\auslass}{(\ldots\hspace{-0.2em})}

\newcommand{\trans}[1]{{#1}^\mathsf{T}}

\newcommand{\bigskipno}{\bigskip \noindent}
\newcommand{\medskipno}{\medskip \noindent}
\newcommand{\smallskipno}{\smallskip \noindent}

\newcommand{\hypref}[1]{(A\ref{#1})}

\newcounter{hypcount}
\newenvironment{hypenv}{\renewcommand{\labelenumi}{(A\arabic{enumi})} \begin{enumerate}\setcounter{enumi}{\value{hypcount}}}
{\setcounter{hypcount}{\value{enumi}}\end{enumerate}}

\theoremstyle{plain}
\newtheorem{thrm}{Theorem}
\newtheorem{prop}{Proposition}
\newtheorem{lemma}{Lemma}
\theoremstyle{definition}
\newtheorem{defn}{Definition}
\theoremstyle{remark}
\newtheorem{case}{Case}


\newenvironment{enumrm}{\begin{enumerate}\renewcommand{\labelenumi}{\textup{(\roman{enumi})}}}{\end{enumerate}}

\renewcommand{\phi}{\varphi}



\title{Discretisation of stochastic control problems for continuous time dynamics with delay\thanks{Financial support by the DFG-Sonderforschungsbereich 649 {\it Economic Risk} is gratefully acknowledged.}}

\author{\sc Markus Fischer \\Weierstra\ss-Institut f\"ur Angewandte\\ Analysis und Stochastik (WIAS)\\ Mohrenstr.\ 39\\ 10117 Berlin\\ Germany \and \sc Markus Rei\ss \\Institut f\"ur Angewandte Mathematik\\ Universit\"at Heidelberg\\ Im Neuenheimer Feld 294\\ 69120 Heidelberg \\ Germany}

\date{February 17, 2006}

\maketitle

\begin{abstract}
As a main step in the numerical solution of control problems in continuous time, the controlled process is approximated by sequences of controlled Markov chains, thus discretising time and space. A new feature in this context is to allow for delay in the dynamics. The existence of an optimal strategy with respect to the cost functional can be guaranteed in the class of relaxed controls. Weak convergence of the approximating extended Markov chains to the original process together with convergence of the associated optimal strategies is established.
\end{abstract}

\section{Introduction} \label{SectIntro}

A general strategy for rendering control problems in continuous time
accessible to numerical computation is the following: Taking as a
starting point the original dynamics, construct a family of control
problems in discrete time with discrete state space and discretised
cost functional. Standard numerical schemes can be applied to find
an optimal control and to calculate the minimal costs for each of
the discrete control problems. The important point to establish is
then whether the discrete optimal controls and minimal costs
converge to the continuous-time limit as the mesh size of the
discretisation tends to zero. If that is the case, then the discrete
control problems are a valid approximation to the original problem.

Approximation schemes for non-delay stochastic control problems in
continuous time implementing the general strategy just outlined are
well established, see \citet{kushnerdupuis01}. The method yields
convergence results under very general conditions.
In the non-delay case \citet{krylov00} derived rates of convergence for those schemes by exploiting fine analytical properties of the associated Bellman equations.

The dynamics of the control problem we are interested in are
described by a stochastic delay differential equation (SDDE). Thus,
the future evolution of the dynamics may depend not only on the
present state, but also on the past evolution. For an exposition of
the general theory of SDDEs see \citet{mohammed84} or \citet{mao97}.
The development of numerical methods for SDDEs has attracted much
attention recently, see \citet{buckwar00}, \citet{huetal04} and the
references therein. In \citet{calzolarietal05} a rate of convergence for a segmentwise Euler scheme is obtained (Proposition 4.2 ibid.) and is used in a non-linear filtering problem for approximating the state process, which is given by an SDDE. Numerical procedures
for deterministic control with delayed dynamics have already been
used in applications, see \citet{boucekkineetal05} for the analysis
of an economic growth model. The algorithm proposed there is based
on the discretisation method studied here, but no formal proof of
convergence is given.

The mathematical analysis of stochastic control problems with time
delay in the state equation has been the object of recent works, see
e.\,g.\ \citet{elsanosietal00} for certain explicitly available
solutions, \citet{oeksendalsulem01} for the derivation of a maximum
principle and \citet{larssen02} for the dynamic programming
approach. Although one can invoke the dynamic programming principle
to derive a Hamilton-Jacobi-Bellman equation for the value function,
such an equation will in general be a non-linear partial
differential equation on a functional state space. The analytical
methods for the non-delay case do not simply carry over to this
infinite-dimensional setting. Another approach to treat stochastic
control problems with delay is based on representing the state
equation as an evolution equation in Hilbert space, see
\citet{bensoussanetal92}.

The class of control problems is specified in
Section~\ref{SectProblem}. In Section~\ref{SectExistence} we prove
the existence of optimal strategies for those problems in the class
of relaxed controls. Section~\ref{SectChainApproximation} introduces
the approximating processes and provides a tightness result.
Finally, in Section~\ref{SectCosts} the discrete control problems
are defined and the convergence of the minimal costs and optimal
strategies is shown.

\section{The control problem} \label{SectProblem}

We consider the control of a dynamical system given by a one-dimensional stochastic delay differential equation (SDDE) driven by a Wiener process. Both drift and diffusion coefficient may depend on the solution's history a certain amount of time into the past. Let $r > 0$ denote the \emph{delay length}, i.\,e.\ the maximal length of dependence on the past. For simplicity, we restrict attention to the case, where only the drift term can be directly controlled.

Typically, the solution process of an SDDE does not enjoy the Markov property, while the segment process associated with that solution does. For a real-valued {\it c\`adl\`ag} function (i.\,e., right-continuous function with left-hand limits) $\psi$ living on the time interval $[-r,\infty)$ the \emph{segment} at time $t \in [0,\infty)$ is defined to be the function $\psi_{t}: [-r,0] \rightarrow \mathbb{R}$ given by $\psi_{t}(s):=\; \psi(t\!+\!s)$. Thus, the segment process $(X_{t})_{t\geq0}$ associated with a real-valued c\`adl\`ag process $(X(t))_{t\geq-r}$ takes its values in $D_{0}\!:= D([-r,0])$, the space of all real-valued c\`adl\`ag functions on the interval $[-r,0]$. There are two natural topologies on $D_{0}$. The first is the one induced by the supremum norm. The second is the \emph{Skorohod topology} of c\`adl\`ag convergence \citep[e.\,g.][]{billingsley99}. The main difference between the Skorohod and the uniform topology lies in the different evaluation of convergence of functions with jumps, which appear naturally as initial segments and discretised processes. For continuous functions both topologies coincide. Similar statements hold for $D_{\infty}\!:= D([-r,\infty))$ and $\tilde{D}_{\infty}\!:= D([0,\infty))$, the spaces of all real-valued c\`adl\`ag functions on the intervals $[-r,\infty)$ and $[0,\infty)$, respectively. The spaces $D_{\infty}$ and $\tilde{D}_{\infty}$ will always be supposed to carry the Skorohod topology, while $D_{0}$ will canonically be equipped with the uniform topology.

Let $(\Gamma,d_{\Gamma})$ be a compact metric space, the space of \emph{control actions}. Denote by $b$ the drift coefficient of the controlled dynamics, and by $\sigma$ the diffusion coefficient. Let $(W(t))_{t\geq0}$ be a one-dimensional standard Wiener process on a filtered probability space $(\Omega,\mathcal{F},(\mathcal{F}_{t})_{t\geq0},\Prb)$ satisfying the usual conditions, and let $(u(t))_{t\geq0}$ be a \emph{control process}, i.\,e.\ an $(\mathcal{F}_{t})$-adapted measurable process with values in $\Gamma$. Consider the controlled SDDE
\begin{equation} \label{EqControlSDDE}
    dX(t) \;=\; b\bigl(X_{t},u(t)\bigr)\,dt \;+\; \sigma(X_{t})\,dW(t), \qquad t \geq 0.
\end{equation}
The control process $u(.)$ together with its stochastic basis including the Wiener process is called an \emph{admissible control} if, for every deterministic initial condition $\phi \in D_{0}$, Eq.~\eqref{EqControlSDDE} has a unique solution which is also weakly unique. Write $\mathcal{U}_{ad}$ for the set of admissible controls of Eq.~\eqref{EqControlSDDE}. The stochastic basis coming with an admissible control will often be omitted in the notation.

A \emph{solution} in the sense used here is an adapted c\`adl\`ag process defined on the stochastic basis of the control process such that the integral version of Eq.~\eqref{EqControlSDDE} is satisfied. Given a control process together with a standard Wiener process, a solution to Eq.~\eqref{EqControlSDDE} is \emph{unique} if it is indistinguishable from any other solution almost surely satisfying the same initial condition. A solution is \emph{weakly unique} if it has the same law as any other solution with the same initial distribution and satisfying Eq.~\eqref{EqControlSDDE} for a control process on a possibly different stochastic basis so that the joint distributions of control and driving Wiener process are the same for both solutions. Let us specify the regularity assumptions to be imposed on the coefficients $b$ and $\sigma$:
\begin{hypenv}
    \item \label{HypCadlag} C\`adl\`ag functionals: the mappings
    \begin{align*}
        & (\psi,\gamma) \mapsto \bigl[t \mapsto b(\psi_{t},\gamma),\; t \geq 0 \bigr], & & \psi \mapsto \bigl[t \mapsto \sigma(\psi_{t}),\; t \geq 0 \bigr]
    \end{align*}
    define measurable functionals $D_{\infty} \times \Gamma \rightarrow \tilde{D}_{\infty}$ and $D_{\infty} \rightarrow \tilde{D}_{\infty}$, respectively, where $D_{\infty}$, $\tilde{D}_{\infty}$ are equipped with their Borel $\sigma$-algebras.
    \item \label{HypDrift} Continuity of the drift coefficient: there is an at most countable subset of $[-r,0]$, denoted by $I_{ev}$, such that for every $t \geq 0$ the function defined by
    \[ D_{\infty}\times\Gamma \ni (\psi,\gamma) \mapsto b(\psi_{t},\gamma) \]
    is continuous on $D_{ev}(t)\times\Gamma$ uniformly in $\gamma \in \Gamma$, where
    \[ D_{ev}(t):=\; \{\psi \in D_{\infty} \;|\; \psi \;\text{is continuous at}\; t+s\; \text{for all}\; s \in I_{ev}\}. \]
    \item \label{HypGrowth} Global boundedness: $|b|$, $|\sigma|$ are bounded by a constant $K > 0$.
    \item \label{HypLipschitz} Uniform Lipschitz condition: There is a constant $K_{L}>0$ such that for all $\phi, \tilde{\phi} \in D_{0}$, all $\gamma \in \Gamma$
    \[ |b(\phi,\gamma)-b(\tilde{\phi},\gamma)| \;+\; |\sigma(\phi)-\sigma(\tilde{\phi})| \;\leq\; K_{L}\cdot\!\sup_{s\in[-r,0]}|\phi(s)-\tilde{\phi}(s)|. \]
    \item \label{HypEllipticity} Ellipticity of the diffusion coefficient: $\sigma(\phi) \geq \sigma_{0}$ for all $\phi \in D_{0}$, where $\sigma_{0}>0$ is a positive constant.
\end{hypenv}
Assumptions \hypref{HypCadlag} and \hypref{HypLipschitz} on the coefficients allow us to invoke Theorem~V.7 in \citet[p.253]{protter03}, which guarantees the existence of a unique solution to Eq.~\eqref{EqControlSDDE} for every piecewise constant control attaining only finitely many different values. The boundedness Assumption \hypref{HypGrowth} poses no limitation except for the initial conditions, because the state evolution will be stopped when the state process leaves a bounded interval. Assumption \hypref{HypDrift} allows us to use ``segmentwise approximations'' of the solution process, see the proof of Proposition~\ref{PropRelaxedCompactness}. The assumptions imposed on the drift coefficient $b$ are satisfied, for example, by
\begin{equation} \label{ExmplDrift}
    b(\phi,\gamma):=\; f\Bigl(\phi(r_{1}),\ldots,\phi(r_{n}), \int_{-r}^{0}\phi(s)w_{1}(s)ds,\ldots,\int_{-r}^{0}\phi(s)w_{m}(s)ds\Bigr)\cdot g(\gamma),
\end{equation}
where $r_{1},\ldots,r_{n} \in [-r,0]$ are fixed, $f$, $g$ are bounded continuous functions and $f$ is Lip\-schitz, and the weight functions $w_{1},\ldots,w_{m}$ lie in $L^1([-r,0])$. Apart from the control term, the diffusion coefficient $\sigma$ may have the same structure as $b$ in \eqref{ExmplDrift}.

We next give an example of a function that could be taken for $\sigma$ if the c\`adl\`ag continuity in Assumption~\hypref{HypCadlag} were missing. In Section~\ref{SectChainApproximation} it will become clear that the corresponding control problem cannot be approximated by a simple discretisation procedure, because the evaluation of $\sigma(\phi)$ for any $\phi \in D_{0}$ depends on the discretisation grid. Let $A_{M}$ be the subset of the interval $[-r,0]$ given by
    \[ A_{M}:=\; \bigl\{ (t-2^{-3M},t] \;\big|\; t = r(\tfrac{n}{2^{M}}-1) \;\text{for some}\; n\in\{1,\ldots,2^{M}\}\bigr\}. \]
Let $A$ be the union of the sets $A_{M}$, $M \in \mathbb{N}$. With positive constants $\sigma_{0}$, $K$, we define a functional $\sigma: D_{0} \rightarrow \mathbb{R}$ by
\begin{equation} \label{ExmplDiffusion}
    \sigma(\phi):=\; \sigma_{0} \;+\; K\wedge\sup \bigl\{|\phi(t)-\phi(t-)| \;\big|\; t \in A \bigr\},
\end{equation}
where $\phi(t-)$ is the left hand limit of $\phi$ at $t \in [-r,0]$. Assumptions \hypref{HypGrowth} and \hypref{HypLipschitz} are clearly satisfied if we choose $\sigma$ according to \eqref{ExmplDiffusion}, but $\sigma$ would not induce a c\`adl\`ag functional $D_{\infty} \rightarrow \tilde{D}_{\infty}$. This can be seen by considering the mapping $[0,\infty) \ni t \mapsto \sigma(\psi_{t})$ for a function $\psi \in D_{\infty}$ which is constant except for a single discontinuity. If we had defined $\sigma$ with the set $A$ being the union of only finitely many sets $A_{M}$, then we would have obtained a c\`adl\`ag functional.

We consider control problems in the weak formulation \citep[cf.][p.\,64]{yongzhou99}. Given an admissible control $u(.)$ and a deterministic initial segment $\phi \in D_{0}$, denote by $X^{\phi,u}$ the unique solution to Eq.~\eqref{EqControlSDDE}. Let $I$ be a compact interval with non-empty interior. Define the stopping time $\tau^{\bar{T}}_{\phi,u}$ of first exit from the interior of $I$ before time $\bar{T}>0$ by
\begin{equation} \label{ExEndTime}
    \tau^{\bar{T}}_{\phi,u} :=\; \inf\{t \geq 0 \;|\; X^{\phi,u}(t) \notin \intr(I)\} \;\wedge\; \bar{T}.
\end{equation}
In order to define the costs, we prescribe a \emph{cost rate} $k\!: \mathbb{R} \times \Gamma \rightarrow [0,\infty)$ and a \emph{boundary cost} $g\!: \mathbb{R} \rightarrow [0,\infty)$ which we take to be (jointly) continuous bounded functions. Let $\beta \geq 0$ denote the exponential \emph{discount rate}. Then define the \emph{cost functional} on $D_{0} \times \mathcal{U}_{ad}$ by
\begin{equation} \label{ExCostFunctional}
    J(\phi,u) :=\; \Mean\left(\int_{0}^{\tau} \exp(-\beta s)\cdot k\bigl(X^{\phi,u}(s),u(s)\bigr)\,ds \;+\; g\bigl(X^{\phi,u}(\tau)\bigr)\right),
\end{equation}
where $\tau = \tau^{\bar{T}}_{\phi,u}$. Our aim is to minimize $J(\phi,.)$. We introduce the \emph{value function}
\begin{equation} \label{ExValueFunction}
    V(\phi) :=\; \inf\{J(\phi,u) \;|\; u \in \mathcal{U}_{ad}\}, \quad \phi \in
D_{0}.
\end{equation}
The control problem now consists in calculating the function $V$ and finding admissible controls that minimize $J$. Such control processes are called \emph{optimal controls} or \emph{optimal strategies}.

\section{Existence of optimal strategies} \label{SectExistence}

In the class $\mathcal{U}_{ad}$ of admissible controls it may happen that there is no optimal control \citep[cf.][p.\,86]{kushnerdupuis01}. A way out is to enlarge the class of controls, allowing for so-called relaxed controls, so that the existence of an optimal (relaxed) control is guaranteed, while the infimum of the costs over the new class coincides with the value function $V$ as given by \eqref{ExValueFunction}.

A \emph{deterministic relaxed control} is a positive measure $\rho$ on $\mathcal{B}(\Gamma\times[0,\infty))$, the Borel $\sigma$-algebra on $\Gamma\times[0,\infty)$, such that
\begin{equation} \label{ExDetermRelaxed}
    \rho(\Gamma\times[0,t]) \;=\; t \quad\text{for all}\; t \geq 0.
\end{equation}
For each $G \in \mathcal{B}(\Gamma)$, the function $t \mapsto \rho(G\times[0,t])$ is absolutely continuous with respect to Lebesgue measure on $[0,\infty)$ by virtue of property \eqref{ExDetermRelaxed}. Denote by $\dot{\rho}(.,G)$ any Lebesgue density of $\rho(G\times[0,.])$. The family of densities $\dot{\rho}(.,G)$, $G \in \mathcal{B}(\Gamma)$, can be chosen in a Borel measurable way such that $\dot{\rho}(t,.)$ is a probability measure on $\mathcal{B}(\Gamma)$ for each $t \geq 0$, and
\begin{equation*}
    \rho(B) \;=\; \int_{0}^{\infty} \int_{\Gamma} \mathbf{1}_{\{(\gamma,t)\in B\}}\, \dot{\rho}(t,d\gamma)\,dt \quad\text{for all}\; B \in \mathcal{B}(\Gamma\times[0,\infty)).
\end{equation*}
Denote by $\mathcal{R}$ the space of deterministic relaxed controls which is equipped with the \emph{weak-compact topology} induced by the following notion of convergence: a sequence
$(\rho_{n})_{n\in\mathbb{N}}$ of relaxed controls converges to $\rho \in \mathcal{R}$ if and only if
\begin{equation*}
    \int\limits_{\Gamma\times[0,\infty)} g(\gamma,t)\, d\rho_{n}(\gamma,t) \;\stackrel{n\to\infty}{\longrightarrow}\; \int\limits_{\Gamma\times[0,\infty)} g(\gamma,t)\, d\rho(\gamma,t) \quad\text{for all}\; g \in \mathbf{C}_{c}(\Gamma\times[0,\infty)),
\end{equation*}
where $\mathbf{C}_{c}(\Gamma\times[0,\infty))$ is the space of all real-valued continuous functions on $\Gamma\times[0,\infty)$ having compact support. Under the weak-compact topology, $\mathcal{R}$ is a (sequentially) compact space.

Suppose $(\rho_{n})_{n\in\mathbb{N}}$ is a convergent sequence in $\mathcal{R}$ with limit $\rho$. Given $T > 0$, let $\rho_{n|T}$ denote the restriction of $\rho_{n}$ to the Borel $\sigma$-algebra on $\Gamma\times[0,T]$, and denote by $\rho_{|T}$ the restriction of
$\rho$ to $\mathcal{B}(\Gamma\times[0,T])$. Then $\rho_{n|T}$, $n\in\mathbb{N}$, $\rho_{|T}$ are all finite measures and $(\rho_{n|T})$ converges weakly to $\rho_{|T}$.

A \emph{relaxed control process} is an $\mathcal{R}$-valued random variable $R$ such that the mapping $\omega \mapsto R(\Gamma\times[0,t])(\omega)$ is $\mathcal{F}_{t}$-measurable for all $t \geq 0$, $G \in \mathcal{B}(\Gamma)$. For a relaxed control process $R$ Eq.~\eqref{EqControlSDDE} takes on the form
\begin{equation} \label{EqRelControlSDDE}
    dX(t) \;=\; \Bigl(\int_{\Gamma} b(X_{t},\gamma)\,\dot{R}(t,d\gamma)\Bigr)dt \;+\; \sigma(X_{t})\,dW(t), \qquad t \geq 0,
\end{equation}
where $(\dot{R}(t,.))_{t\geq0}$ is the family of derivative measures associated with $R$. The family $(\dot{R}(t,.))$ can be constructed in a measurable way \citep[cf.][p.\,52]{kushner90}. A relaxed control process together with its stochastic basis including the Wiener process is called \emph{admissible relaxed control} if, for every deterministic initial condition, Eq.~\eqref{EqRelControlSDDE} has a unique solution which is also weakly unique. Any ordinary control process $u$ can be represented as a relaxed control process by setting
\begin{equation*}
    R(B):=\; \int_{0}^{\infty} \int_{\Gamma} \mathbf{1}_{\{(\gamma,t)\in B\}}\, \delta_{u(t)}(d\gamma)\,dt, \quad B \in \mathcal{B}(\Gamma\times[0,\infty)),
\end{equation*}
where $\delta_{\gamma}$ is the Dirac measure at $\gamma \in \Gamma$. Denote by $\hat{\mathcal{U}}_{ad}$ the set of all admissible relaxed controls. Instead of \eqref{ExCostFunctional} we define a cost functional on $D_{0}\times\hat{\mathcal{U}}_{ad}$ by
\begin{equation} \label{ExRelCostFunctional}
    \hat{J}(\phi,R) :=\; \Mean\left(\int_{0}^{\tau}\int_{\Gamma} \exp(-\beta s)\cdot k\bigl(X^{\phi,R}(s),\gamma\bigr)\,\dot{R}(s,d\gamma)\,ds \;+\; g\bigl(X^{\phi,R}(\tau)\bigr)\right),
\end{equation}
where $X^{\phi,R}$ is the solution to Eq.~\eqref{EqRelControlSDDE} under the relaxed control process $R$ with initial segment $\phi$ and $\tau$ is defined in analogy to \eqref{ExEndTime}. Instead of \eqref{ExValueFunction} as value function we have
\begin{equation} \label{ExRelValueFunction}
    \hat{V}(\phi) :=\; \inf\{\hat{J}(\phi,R) \;|\; R \in \hat{\mathcal{U}}_{ad}\}, \quad \phi \in D_{0}.
\end{equation}
The cost functional $\hat{J}$ depends only on the joint distribution of the solution $X^{\phi,R}$ and the underlying control process $R$, since $\tau$, the time horizon, is a deterministic function of the solution. The distribution of $X^{\phi,R}$, in turn, is determined by the initial condition $\phi$ and the joint distribution of the control process and its accompanying Wiener process. Letting the time horizon vary, we may regard $\hat{J}$ as a function of the law of $(X,R,W,\tau)$, that is, as being defined on a subset of the set of probability measures on $\mathcal{B}(D_{\infty}\times\mathcal{R}\times \tilde{D}_{\infty}\times[0,\infty])$. Notice that the time interval has been compactified. The domain of definition of $\hat{J}$ is determined by the class of admissible relaxed controls for Eq.~\eqref{EqRelControlSDDE}, the definition of the time horizon and the distributions of the initial segments $X_{0}$.

The idea in proving existence of an optimal strategy is to check that $\hat{J}(\phi,.)$ is a (sequentially) lower semi-continuous function defined on a (sequentially) compact set. It then follows from a theorem by Weierstra{\ss} that $\hat{J}(\phi,.)$ attains its minimum at some point of its compact domain \citep[cf.][p.\,65]{yongzhou99}. The following proposition gives the analogue of Theorem 10.1.1 in \citet[pp.\,271-275]{kushnerdupuis01} for our setting. We present the proof in detail, because the identification of the limit process is different from the classical case.

\begin{prop} \label{PropRelaxedCompactness}
Assume \hypref{HypCadlag}\,--\,\hypref{HypLipschitz}. Let $((R^{M},W^{M}))_{M\in\mathbb{N}}$ be any sequence of admissible relaxed controls for Eq.~\eqref{EqRelControlSDDE}, where $(R^{M},W^{M})$ is defined on the filtered probability space $(\Omega_{M},\mathcal{F}^{M},(\mathcal{F}^{M}_{t}),\Prb_{M})$. Let $X^{M}$ be a solution to Eq.~\eqref{EqRelControlSDDE} under control $(R^{M},W^{M})$ with deterministic initial condition $\phi^{M} \in D_{0}$, and assume that $(\phi^{M})$ tends to $\phi$ uniformly for some $\phi \in D_{0}$. For each $M \in \mathbb{N}$, let $\tau^{M}$ be an $(\mathcal{F}^{M}_{t})$-stopping time. Then $((X^{M},R^{M},W^{M},\tau^{M}))_{M\in\mathbb{N}}$ is tight.

Denote by $(X,R,W,\tau)$ a limit point of the sequence $((X^{M},R^{M},W^{M},\tau^{M}))_{M\in\mathbb{N}}$. Define a filtration by $\mathcal{F}_{t}\!:= \sigma(X(s),R(s),W(s),\tau\mathbf{1}_{\{\tau\leq t\}},\, s \leq t)$, $t \geq 0$. Then $W(.)$ is an $(\mathcal{F}_{t})$-adapted Wiener process, $\tau$ is an $(\mathcal{F}_{t})$-stopping time, $(R,W)$ is an admissible relaxed control, and $X$ is a solution to Eq.~\eqref{EqRelControlSDDE} under $(R,W)$ with initial condition $\phi$.
\end{prop}

\begin{proof}
Tightness of $(X^{M})$ follows from the Aldous criterion \citep[cf.][pp.\,176-179]{billingsley99}: given $M \in \mathbb{N}$, any bounded $(\mathcal{F}^{M}_{t})$-stopping time $\nu$ and $\delta > 0$ we have
\begin{equation*}
    \Mean_{M}\bigl(\bigl|X^{M}(\nu+\delta)-X^{M}(\nu)\bigr|^{2} \;\big|\; \mathcal{F}^{M}_{\nu}\bigr) \quad\leq\quad 2K^{2}\delta(\delta+1)
\end{equation*}
as a consequence of Assumption \hypref{HypGrowth} and the It\^o isometry. Notice that $X^{M}(0)$ tends to $X(0)$ as $M$ goes to infinity by hypothesis. The sequences $(R^{M})$ and $(\tau^{M})$ are tight, because the value spaces $\mathcal{R}$ and $[0,\infty]$, respectively, are compact. The sequence $(W^{M})$ is tight, since all $W^{M}$ induce the same measure. Finally, componentwise tightness implies tightness of the product \citep[cf.][p.\,65]{billingsley99}.

By abuse of notation, we do not distinguish between the convergent subsequence and the original sequence and assume that $((X^{M},R^{M},W^{M},\tau^{M}))$ converges weakly to $(X,R,W,\tau)$. The random time $\tau$ is an $(\mathcal{F}_{t})$-stopping time by construction of the filtration. Likewise, $R$ is $(\mathcal{F}_{t})$-adapted by construction, and it is indeed a relaxed control process, because $R(\Gamma\times[0,t]) = t$, $t \geq 0$, $\Prb$-almost surely by weak convergence of the relaxed control processes $(R^{M})$ to $R$. The process $W$ has Wiener distribution and continuous paths with probability one, being the limit of standard Wiener processes. To check that $W$ is an $(\mathcal{F}_{t})$-Wiener process, we use the martingale problem characterization of Brownian motion. To this end, for $g \in \mathbf{C}_{c}(\Gamma\times[0,\infty))$, $\rho \in \mathcal{R}$ define the pairing
\begin{equation*}
    (g,\rho)(t):=\; \int_{\Gamma\times[0,t]} g(\gamma,s)\,d\rho(\gamma,s), \quad t \geq 0.
\end{equation*}
Notice that real-valued continuous functions on $\mathcal{R}$ can be approximated by functions of the form
\begin{equation*}
    \mathcal{R} \ni \rho \mapsto \tilde{H}\bigl((g_{j},\rho)(t_{i}),\, (i,j)\in\mathbb{N}_{p}\times\mathbb{N}_{q}\bigr)\in\mathbb{R},
\end{equation*}
where $p$, $q$ are natural numbers, $\{t_{i} \;|\; i\in\mathbb{N}_{p}\} \subset [0,\infty)$, and $\tilde{H}$, $g_{j}$, $j \in \mathbb{N}_{q}$, are suitable continuous functions with compact support and $\mathbb{N}_{N}:=\{1,\ldots,N\}$ for any $N\in\mathbb{N}$. Let $t \geq 0$, $t_{1},\ldots,t_{p} \in [0,t]$, $h \geq 0$, $g_{1},\ldots,g_{q}$ be functions in $\mathbf{C}_{c}(\Gamma\times[0,\infty))$, and $H$ be a continuous function of $2p+p\!\cdot\!q+1$ arguments with compact support. Since $W^{M}$ is an $(\mathcal{F}^{M}_{t})$-Wiener process for each $M \in \mathbb{N}$, we have for all $f \in \mathbf{C}^{2}_{c}(\mathbb{R})$
\begin{equation*} \begin{split}
    \Mean_{M}\Bigl(&H\bigl(X^{M}(t_{i}),(g_{j},R^{M})(t_{i}),W^{M}(t_{i}), \tau^{M}\mathbf{1}_{\{\tau^{M}\leq t\}},\, (i,j)\in\mathbb{N}_{p}\times\mathbb{N}_{q}\bigr) \\
    &\cdot\Bigl(f\bigl(W^{M}(t+h)\bigr) - f\bigl(W^{M}(t)\bigr) - \frac{1}{2} \int\limits_{t}^{t+h} \frac{\partial^{2}f}{\partial x^{2}}\bigl(W^{M}(s)\bigr)ds \Bigr)\Bigr) \quad=\quad 0.
    \end{split}
\end{equation*}
By the weak convergence of $((X^{M},R^{M},W^{M},\tau^{M}))_{M\in\mathbb{N}}$ to $(X,W,R,\tau)$ we see that
\begin{equation*} \begin{split}
    \Mean\Bigl(&H\bigl(X(t_{i}),(g_{j},R)(t_{i}),W(t_{i}), \tau\mathbf{1}_{\{\tau\leq t\}},\, (i,j)\in\mathbb{N}_{p}\times\mathbb{N}_{q}\bigr) \\
    &\cdot\Bigl(f\bigl(W(t+h)\bigr) - f\bigl(W(t)\bigr) - \frac{1}{2} \int\limits_{t}^{t+h} \frac{\partial^{2}f}{\partial x^{2}}\bigl(W(s)\bigr)ds \Bigr)\Bigr) \quad=\quad 0
\end{split}
\end{equation*}
for all $f \in \mathbf{C}^{2}_{c}(\mathbb{R})$. As $H$, $p$, $q$, $t_{i}$, $g_{j}$ vary over all possibilities, the corresponding random variables $H(X(t_{i}),(g_{j},R)(t_{i}),W(t_{i}), \tau\mathbf{1}_{\{\tau\leq t\}},\, (i,j)\in\mathbb{N}_{p}\times\mathbb{N}_{q})$ induce the $\sigma$-algebra $\mathcal{F}_{t}$. Since $t \geq 0$, $h \geq 0$ were arbitrary, it follows that
\begin{equation*}
	f\bigl(W(t)\bigr) - f\bigl(W(0)\bigr) - \frac{1}{2} \int\limits_{0}^{t} \frac{\partial^{2}f}{\partial x^{2}}\bigl(W(s)\bigr)ds, \quad t \geq 0,
\end{equation*}
is an $(\mathcal{F}_{t})$-martingale for every $f \in \mathbf{C}^{2}_{c}(\mathbb{R})$. Consequently, $W$ is an $(\mathcal{F}_{t})$-Wiener process.

It remains to show that $X$ solves Eq.~\eqref{EqRelControlSDDE} under control $(R,W)$ with initial condition $\phi$. Notice that $X$ has continuous paths on $[0,\infty)$ $\Prb$-almost surely, because the process $(X(t))_{t\geq0}$ is the weak limit in $\tilde{D}_{\infty}$ of continuous processes. Fix $T > 0$. We have to check that $\Prb$-almost surely
\begin{equation*}
X(t)=\phi(0) \;+\; \int_{0}^{t} \int_{\Gamma} b(X_{s},\gamma)\,\dot{R}(s,d\gamma)\,ds
+\; \int_{0}^{t} \sigma(X_{s})\,dW(s) \text{ for all}\; t \in [0,T].
\end{equation*}
By virtue of the Skorohod representation theorem \citep[cf.][p.\,70]{billingsley99} we may assume that the processes $(X^{M},R^{M},W^{M})$, $M \in \mathbb{N}$, are all defined on the same probability space $(\Omega,\mathcal{F},\Prb)$ as $(X,R,W)$ and that convergence of $((X^{M},R^{M},W^{M}))$ to $(X,R,W)$ is $\Prb$-almost sure. Since $X$, $W$ have continuous paths on $[0,T]$ and $(\phi^{M})$ converges to $\phi$ in the uniform topology, one finds $\tilde{\Omega} \in \mathcal{F}$ with $\Prb(\tilde{\Omega}) = 1$ such that for all $\omega \in \tilde{\Omega}$
\begin{align*}
    & \sup_{t\in[-r,T]}\bigl|X^{M}(t)(\omega)-X(t)(\omega)\bigr| \;\stackrel{M\to\infty}{\longrightarrow}\; 0,\! & & \sup_{t\in[-r,T]}\bigl|W^{M}(t)(\omega)-W(t)(\omega)\bigr| \;\stackrel{M\to\infty}{\longrightarrow}\; 0,\! &
\end{align*}
and also $R^{M}(\omega) \to R(\omega)$ in $\mathcal{R}$. Let $\omega \in \tilde{\Omega}$. We first show that
\begin{equation*}
    \int_{0}^{t} \int_{\Gamma} b\bigl(X^{M}_{s}(\omega),\gamma\bigr)\, \dot{R}^{M}(s,d\gamma)(\omega)\,ds \;\stackrel{M\to\infty}{\longrightarrow}\; \int_{0}^{t} \int_{\Gamma} b\bigl(X_{s}(\omega),\gamma\bigr)\,\dot{R}(s,d\gamma)(\omega)\,ds
\end{equation*}
uniformly in $t \in [0,T]$. As a consequence of Assumption~\hypref{HypLipschitz}, the uniform convergence of the trajectories on $[-r,T]$ and property \eqref{ExDetermRelaxed} of the relaxed controls, we have
\begin{equation*}
    \int_{\Gamma\times[0,T]} \bigl|b\bigl(X^{M}_{s}(\omega),\gamma\bigr)-b\bigl(X_{s}(\omega),\gamma\bigr)\bigr|\, dR^{M}(\gamma,s)(\omega) \;\stackrel{M\to\infty}{\rightarrow}\; 0.
\end{equation*}
By Assumption~\hypref{HypDrift}, we find a countable set $A_{\omega} \subset [0,T]$ such that the mapping $(\gamma,s) \mapsto b(X_{s}(\omega),\gamma)$ is continuous in all $(\gamma,s) \in \Gamma\times([0,T]\setminus A_{\omega})$. Since $A_{\omega}$ is countable we have $R(\omega)(\Gamma\times A_{\omega}) = 0$. Hence, by the generalized mapping theorem
\citep[cf.][p.\,21]{billingsley99}, we obtain for each $t \in [0,T]$
\begin{equation*}
    \int_{\Gamma\times[0,t]} b\bigl(X_{s}(\omega),\gamma\bigr)\, dR^{M}(\gamma,s)(\omega) \;\stackrel{M\to\infty}{\rightarrow}\; \int_{\Gamma\times[0,t]} b\bigl(X_{s}(\omega),\gamma\bigr)\, dR(\gamma,s)(\omega).
\end{equation*}
The convergence is again uniform in $t \in [0,T]$, as $b$ is bounded and $R^{M}$, $M \in \mathbb{N}$, $R$ are all positive measures with mass $T$ on $\Gamma\times[0,T]$. Define c\`adl\`ag processes $C^{M}$, $M \in \mathbb{N}$, on $[0,\infty)$ by
\begin{equation*}
    C^{M}(t):=\; \phi^{M}(0) \;+\; \int_{\Gamma\times[0,t]} b(X^{M}_{s},\gamma)\, dR^{M}(\gamma,s),\quad t \geq 0,
\end{equation*}
and define $C$ in analogy to $C^{M}$ with $\phi$, $R$, $X$ in place of $\phi^{M}$, $R^{M}$, $X^{M}$, respectively. From the above, we know that $C^{M}(t)\to C(t)$ holds uniformly over $t\in[0,T]$ for any $T > 0$ with probability one. Define operators $F^{M}\!:\tilde{D}_{\infty} \rightarrow \tilde{D}_{\infty}$, $M\in\mathbb{N}$, mapping c\`adl\`ag processes to c\`adl\`ag processes by
\begin{equation*}
    F^{M}(Y)(t)(\omega):=\; \sigma\left([-r,0] \ni s \mapsto \begin{cases} Y(t\!+\!s)(\omega) &\text{if}\; t\!+\!s \geq 0, \\ \phi^{M}(t\!+\!s) &\text{else}\end{cases} \right),\quad t \geq 0,\; \omega \in \Omega,
\end{equation*}
and define $F$ in the same way as $F^{M}$ with $\phi^{M}$ replaced by $\phi$. Observe that $X^{M}$ solves
\begin{equation*}
    X^{M}(t) \;=\; C^{M}(t) \;+\; \int_{0}^{t} F^{M}(X^{M})(s-)\,dW^{M}(s), \quad t \geq 0.
\end{equation*}
Denote by $(\hat{X}(t))_{t\geq0}$ the unique solution to
\begin{equation*}
    \hat{X}(t) \;=\; C(t) \;+\; \int_{0}^{t} F(\hat{X})(s-)\,dW(s), \quad t \geq 0,
\end{equation*}
and set $\hat{X}(t)\!:=\phi(t)$ for $t \in [-r,0)$. Assumption~\hypref{HypLipschitz} and the uniform convergence of $(\phi^{M})$ to $\phi$ imply that $F^{M}(\hat{X})$ converges to $F(\hat{X})$ uniformly on compacts in probability (convergence \emph{in ucp}). Theorem~V.15 in \citet[p.\,265]{protter03} yields that $(X^{M})$ converges to $\hat{X}$ in ucp, that is
\begin{equation*}
    \sup_{t\in[0,T]}\bigl|X^{M}(t)-\hat{X}(t)\bigr| \;\stackrel{M\to\infty}{\longrightarrow}\; 0 \quad\text{in probability $\Prb$ for any}\; T > 0.
\end{equation*}
Therefore, $X$ is indistinguishable from $\hat{X}$. By definition of $C$ and $F$, this implies that $\hat{X}$ solves Eq.~\eqref{EqRelControlSDDE} under control $(R,W)$ with initial condition $\phi$, and so does $X$.
\end{proof}

If the time horizon were deterministic, then the existence of optimal strategies in the class of relaxed controls would be clear. Given an initial condition $\phi \in D_{0}$, one would select a sequence $((R^{M},W^{M}))_{M\in\mathbb{N}}$ such that $(\hat{J}(\phi,R^{M}))$ converges to its infimum. By Proposition~\ref{PropRelaxedCompactness}, a suitable subsequence of $((R^{M},W^{M}))$ and the associated solution processes would converge weakly to $(R,W)$ and the associated solution to Eq.~\eqref{EqRelControlSDDE}. Taking into account \eqref{ExRelCostFunctional}, the definition of the costs, this in turn would imply that $\hat{J}(\phi,.)$ attains its minimum value at $R$ or, more precisely, $(X,R,W)$.

A similar argument is still valid, if the time horizon depends continuously on the paths with probability one under every possible solution. That is to say, the mapping
\begin{align} \label{ExExitTime}
    & \hat{\tau}:\; D_{\infty} \rightarrow [0,\infty], & & \hat{\tau}(\psi) :=\; \inf\{t \geq 0 \;|\; \psi(t) \notin \intr(I)\} \;\wedge\; \bar{T}, &
\end{align}
is Skorohod continuous with probability one under the measure induced by any solution $X^{\phi,R}$, $R$ any relaxed control. This is indeed the case if the diffusion coefficient $\sigma$ is bounded away from zero as required by Assumption~\hypref{HypEllipticity}, cf.\ \citet[pp.\,277-281]{kushnerdupuis01}.

By introducing relaxed controls, we have enlarged the class of possible strategies. The infimum of the costs, however, remains the same for the new class. This is a consequence of the fact that stochastic relaxed controls can be arbitrarily well approximated by piecewise constant ordinary stochastic controls which attain only a finite number of different control values. A proof of this assertion is given in \citet[pp.\,59-60]{kushner90} in case the time horizon
is finite, and extended to the case of control up to an exit time in \citet[pp.\,282-286]{kushnerdupuis01}. Notice that nothing hinges on the presence or absence of delay in the controlled dynamics. Let us summarize our findings.

\begin{thrm} \label{ThOptimalControl}
Assume \hypref{HypCadlag}\,--\,\hypref{HypEllipticity}. Given any deterministic initial condition $\phi \in D_{0}$, the relaxed control problem determined by \eqref{EqRelControlSDDE} and \eqref{ExRelCostFunctional} possesses an optimal strategy, and the minimal costs are the same as for the original control problem.
\end{thrm}

\section{Approximating chains} \label{SectChainApproximation}

In order to construct finite-dimensional approximations to our
control problem, we discretise time and state space.
Denote by $h>0$ the mesh size of an equidistant time discretisation starting at
zero. Let $S_{h}\!:= \sqrt{h}\mathbb{Z}$ be the corresponding state
space, and set $I_{h}\!:= I \cap S_{h}$. Notice that $S_{h}$ is
countable and $I_{h}$ is finite. Let $\Lambda_{h}\!: \mathbb{R}
\rightarrow S_{h}$ be a round-off function. We will simplify things
even further by considering only mesh sizes $h = \frac{r}{M}$ for
some $M \in \mathbb{N}$, where $r$ is the delay length. The number
$M$ will be referred to as \emph{discretisation degree}.

The admissible controls for the finite-dimensional control problems correspond to piecewise constant processes in continuous time. A time-discrete process $u=(u(n))_{n\in\mathbb{N}_{0}}$ on a stochastic basis $(\Omega,\mathcal{F},(\mathcal{F}_{t}),\Prb)$ with values in $\Gamma$ is a \emph{discrete admissible control of degree $M$} if $u$ takes on only finitely many different values in $\Gamma$ and $u(n)$ is $\mathcal{F}_{nh}$-measurable for all $n \in \mathbb{N}_{0}$. Denote by $(\bar{u}(t))_{t\geq0}$ the piecewise constant c\`adl\`ag interpolation to $u$ on the time grid. We call a time-discrete process $(\xi(n))_{n\in\{-M,\ldots,0\}\cup\mathbb{N}}$ a \emph{discrete chain of degree $M$} if $(\xi(n))$ takes its values in $S_{h}$ and $\xi(n)$ is $\mathcal{F}_{nh}$-measurable for all $n \in \mathbb{N}_{0}$. In analogy to $\bar{u}$, write $(\bar{\xi}(t))_{t\geq-r}$ for the c\`adl\`ag interpolation to the discrete chain $(\xi(n))_{n\in\{-M,\ldots,0\}\cup\mathbb{N}}$. We denote by $\bar{\xi}_{t}$ the $D_{0}$-valued segment of $\bar{\xi}(.)$ at time $t \geq 0$.

Let $\phi \in D_{0}$ be a deterministic initial condition, and suppose we are given a \emph{sequence of discrete admissible controls} $(u^{M})_{M\in\mathbb{N}}$, that is $u^{M}$ is a discrete admissible control of degree $M$ on a stochastic basis $(\Omega_{M},\mathcal{F}^{M},(\mathcal{F}^{M}_{t}),\Prb_{M})$ for each $M \in \mathbb{N}$. In addition, suppose that the sequence $(\bar{u}^{M})$ of interpolated discrete controls converges weakly to some relaxed control $R$. We are then looking for a sequence approximating the solution $X$ of Eq.~\eqref{EqRelControlSDDE} under control $(R,W)$ with initial condition $\phi$, where the Wiener process $W$ has to be constructed from the approximating sequence.

Given $M$-step or \emph{extended Markov transition functions} $p^{M}\!: S_{h}^{M+1}\times \Gamma \times S_{h} \rightarrow [0,1]$, $M \in \mathbb{N}$, we define a \emph{sequence of approximating chains} associated with $\phi$ and $(u^{M})$ as a family $(\xi^{M})_{M\in\mathbb{N}}$ of processes such that $\xi^{M}$ is a discrete chain of degree $M$ defined on the same stochastic basis as $u^{M}$, provided the following conditions are fulfilled for $h = h_M\!:= \frac{r}{M}$ tending to zero:
\begin{enumrm}
    \item Initial condition: $\xi^{M}(n) = \Lambda_{h}(\phi(nh))$ for all $n \in \{-M,\ldots,0\}$.
    \item Extended Markov property: for all $n \in \mathbb{N}_{0}$, all $x \in S_{h}$
    \begin{equation*}
        \Prb_{M}\bigl(\xi^{M}(n\!+\!1) = x \;\big|\; \mathcal{F}^{M}_{nh}\bigr) \quad=\quad p^{M}\bigl(\xi^{M}(n\!-\!M),\ldots,\xi^{M}(n),u^{M}(n),x\bigr).
    \end{equation*}
    \item Local consistency with the drift coefficient:
    \begin{equation*} \begin{split}
        \mu_{\xi^{M}}(n)\;:=\quad &\Mean_{M}\bigl(\xi^{M}(n\!+\!1) - \xi^{M}(n) \;\big|\; \mathcal{F}^{M}_{nh}\bigr) \\[1ex]
        =\quad & h\cdot b\bigl(\bar{\xi}^{M}_{nh},u^{M}(n)\bigr) + o(h) \quad=:\; h\cdot b_{h}\bigl(\bar{\xi}^{M}_{nh},u^{M}(n)\bigr).
        \end{split}
    \end{equation*}
    \item Local consistency with the diffusion coefficient:
    \begin{equation*}
        \Mean_{M}\bigl(\bigl(\xi^{M}(n\!+\!1) - \xi^{M}(n) - \mu_{\xi^{M}}(n)\bigr)^{2}\big|\; \mathcal{F}^{M}_{nh}\bigr) \;=\;h\cdot \sigma^{2}(\bar{\xi}^{M}_{nh}) + o(h) \;=: h\cdot \sigma^{2}_{h}(\bar{\xi}^{M}_{nh}).
\end{equation*}
\item Jump heights: there is a positive number $\tilde{N}$, independent of $M$, such that
    \[\sup_{n}|\xi^{M}(n+1)-\xi^{M}(n)|\leq\tilde{N}\sqrt{h_{M}}.\]
\end{enumrm}
It is straightforward, under Assumptions \hypref{HypGrowth} and \hypref{HypEllipticity}, to construct a sequence of extended Markov transition functions such that the jump height and the local consistency conditions can be fulfilled. Assuming that the bounding constant $K$ from \hypref{HypGrowth} is a natural number, we may define the functions $p^{M}$ for all $M \in \mathbb{N}$ big enough by, for example,
\begin{equation*}
	p^{M}(Z(-M),\ldots,Z(0),\gamma,x):=\; \begin{cases} \frac{1}{2K^{2}}\sigma(\bar{Z}) + \frac{\sqrt{h}}{2K}b(\bar{Z},\gamma), &\text{if}\; x = Z(0) + K\sqrt{h}, \\
		\frac{1}{2K^{2}}\sigma(\bar{Z}) - \frac{\sqrt{h}}{2K}b(\bar{Z},\gamma), &\text{if}\; x = Z(0)-K\sqrt{h}, \\
		1 - \frac{1}{K^{2}}\sigma(\bar{Z}) &\text{if}\; x = Z(0) \\
		0 &\text{else,}
	\end{cases}
\end{equation*}
where $h = h_{M}$, $Z = (Z(-M),\ldots,Z(0)) \in S_{h}^{M+1}$, $\gamma \in \Gamma$, $x \in S_{h}$, and $\bar{Z} \in D_{0}$ is the piecewise constant interpolation associated with $Z$. The family $(p^{M})$ as just defined, in turn, is all we need in order to construct a sequence of approximating chains associated with any given $\phi$, $(u^{M})$.

We will represent the interpolation $\bar{\xi}^{M}$ as a solution to an equation corresponding to Eq.~\eqref{EqControlSDDE} with control process $\bar{u}^{M}$ and initial condition $\phi^{M}$, where $\phi^{M}$ is the piecewise constant $S_{h}$-valued c\`adl\`ag interpolation to $\phi$, that is $\phi^{M} = \bar{\xi}^{M}_{0}$. Define the discrete process $(L^{M}(n))_{n\in\mathbb{N}_{0}}$ by $L^{M}(0)\!:=0$ and
\begin{equation*}
    \xi^{M}(n) \quad=\quad \phi^{M}(0) \;+\; \sum_{i=0}^{n-1} h\cdot b_{h}\bigl(\bar{\xi}^{M}_{ih},u^{M}(i)\bigr) \;+\; L^{M}(n), \qquad n \in \mathbb{N}.
\end{equation*}
Observe that $L^{M}$ is a martingale in discrete time with respect to the filtration $(\mathcal{F}^{M}_{nh})$. Setting
\begin{equation*}
    \varepsilon^{M}_{1}(t):=\; \sum_{i=0}^{\lfloor\frac{t}{h}\rfloor-1} h\cdot b_{h}\bigl(\bar{\xi}^{M}_{ih},\bar{u}^{M}(ih)\bigr) \;-\; \int_{0}^{t} b\bigl(\bar{\xi}^{M}_{s},\bar{u}^{M}(s)\bigr)\,ds, \qquad t \geq 0,
\end{equation*}
the interpolated process $\bar{\xi}^{M}$ can be represented as solution to
\begin{equation*}
    \bar{\xi}^{M}(t) \quad=\quad \phi^{M}(0) \;+\; \int_{0}^{t} b\bigl(\bar{\xi}^{M}_{s},\bar{u}^{M}(s)\bigr)\,ds \;+\; L^{M}(\lfloor\tfrac{t}{h}\rfloor) \;+\; \varepsilon^{M}_{1}(t), \qquad t \geq 0.
\end{equation*}
With $T > 0$, we have for the error term
\begin{equation*} \begin{split}
    \Mean_{M}\Bigl(\sup\nolimits_{t\in[0,T]}\big|\varepsilon^{M}_{1}(t)\big|\Bigr) \quad\leq\quad & \sum_{i=0}^{\lfloor\frac{T}{h}\rfloor-1}\, h\Mean_{M}\Bigl( \bigl|b_{h}\bigl(\bar{\xi}^{M}_{ih},u^{M}(i)\bigr) - b\bigl(\bar{\xi}^{M}_{ih},u^{M}(i)\bigr)\bigr|\Bigr) \;+\; K\cdot h\\[1ex]
+\;& \int_{0}^{h\lfloor\frac{T}{h}\rfloor} \Mean_{M}\Bigl( \bigl|b\bigl(\bar{\xi}^{M}_{h\lfloor\frac{s}{h}\rfloor},\bar{u}^{M}(s)\bigr) - b\bigl(\bar{\xi}^{M}_{s},\bar{u}^{M}(s)\bigr)\bigr|\Bigr)\, ds,
\end{split}
\end{equation*}
which tends to zero as $M$ goes to infinity by Assumptions \hypref{HypDrift}, \hypref{HypGrowth}, dominated convergence and the defining properties of $(\xi^{M})$. Moreover, $|\varepsilon^{M}_{1}(t)|$ is bounded by $2K\!\cdot\!T$ for all $t \in [0,T]$ and all $M$ big enough, whence also
\begin{equation*}
    \Mean_{M}\Bigl(\sup\nolimits_{t\in[0,T]}\bigl|\varepsilon^{M}_{1}(t)\bigr|^{2}\Bigr) \quad\stackrel{M\to\infty}{\longrightarrow}\quad 0.
\end{equation*}
The discrete-time martingale $L^{M}$ can be rewritten as a discrete stochastic integral. Define $(W^{M}(n))_{n\in\mathbb{N}_{0}}$ by setting $W^{M}(0):=0$ and
\begin{equation*}
    W^{M}(n):=\; \sum_{i=0}^{n-1} \frac{1}{\sigma(\bar{\xi}^{M}_{ih})} \bigl(L^{M}(i\!+\!1)-L^{M}(i)\bigr), \qquad n \in \mathbb{N}.
\end{equation*}
Using the piecewise constant interpolation $\bar{W}^{M}$ of $W^{M}$, the process $\bar{\xi}^{M}$ can be expressed as the solution to
\begin{equation} \label{EqApproximateSDDE}
    \bar{\xi}^{M}(t) \;=\; \phi^{M}(0) \;+\; \int_{0}^{t} b\bigl(\bar{\xi}^{M}_{s},\bar{u}^{M}(s)\bigr)\,ds \;+\; \int_{0}^{t} \sigma\Bigl(\bar{\xi}^{M}_{h\lfloor\frac{s-}{h}\rfloor}\Bigr)\, d\bar{W}^{M}(s) \;+\; \varepsilon^{M}_{2}(t), \quad t \geq 0,
\end{equation}
where the error terms $(\varepsilon^{M}_{2})$ converge to zero as $(\varepsilon^{M}_{1})$ before.

We are now prepared for the convergence result, which should be compared to Theorem 10.4.1 in \citet[p.\,290]{kushnerdupuis01}. The proof is similar to that of Proposition~\ref{PropRelaxedCompactness}. We merely point out the main differences.

\begin{prop} \label{PropChainConvergence}
Assume \hypref{HypCadlag}\,--\,\hypref{HypEllipticity}. For each $M
\in \mathbb{N}$, let $\tau^{M}$ be a stopping time with respect to
the filtration generated by
$(\bar{\xi}^{M}(s),\bar{u}^{M}(s),\bar{W}^{M}(s),\, s \leq t)$. Let
$R^{M}$ denote the relaxed control representation of $\bar{u}^{M}$.
Suppose $(\phi^M)$ converges to the initial condition $\phi$
uniformly on $[-r,0]$. Then
$((\bar{\xi}^{M},R^{M},\bar{W}^{M},\tau^{M}))_{M\in\mathbb{N}}$ is
tight.

For a limit point $(X,R,W,\tau)$ set $\mathcal{F}_{t}\!:=
\sigma\bigl(X(s),R(s),W(s),\tau\mathbf{1}_{\{\tau\leq t\}},\; s \leq
t\bigr)$, $t \geq 0$. Then $W$ is an $(\mathcal{F}_{t})$-adapted
Wiener process, $\tau$ is an $(\mathcal{F}_{t})$-stopping time,
$(R,W)$ is an admissible relaxed control, and $X$ is a solution to
Eq.~\eqref{EqRelControlSDDE} under $(R,W)$ with initial
condition $\phi$.
\end{prop}

\begin{proof}
The main differences in the proof are establishing the tightness of $(\bar{W}^{M})$ and the identification of the limit points. We calculate the order of convergence for the discrete-time previsible quadratic variations of $(W^{M})$:
\begin{equation*}
    \langle W^{M}\rangle_{n} \;=\; \sum_{i=0}^{n-1} \Mean\bigl((W^{M}(i\!+\!1)-W^{M}(i))^{2} \;\big|\; \mathcal{F}^{M}_{ih}\bigr) \;=\; n h \;+\; o(h)\sum_{i=0}^{n-1} \frac{1}{\sigma^{2}(\tilde{\xi}^{M}_{ih})}
\end{equation*}
for all $M \in \mathbb{N}$, $n \in \mathbb{N}_{0}$. Taking into account Assumption~\hypref{HypEllipticity} and the definition of the time-continuous processes $\bar{W}^{M}$, we see that $\langle\bar{W}^{M}\rangle$ tends to $\Id_{[0,\infty)}$ in probability uniformly on compact time intervals. By Theorem~VIII.3.11 of \citet[p.\,432]{jacodshiryaev87} we conclude that $(\bar{W}^{M})$ converges weakly in $\tilde{D}_{\infty}$ to a standard Wiener process $W$. That $W$ has independent increments with respect to the filtration $(\mathcal{F}_{t})$ can be seen by considering the first and second conditional moments of the increments of $W^{M}$ for each $M\in\mathbb{N}$ and applying the conditions on local consistency and the jump heights of $(\xi^{M})$.

Suppose $((\bar{\xi}^{M},R^{M},\bar{W}^{M}))$ is weakly convergent with limit point $(X,R,W)$. The remaining different part is the identification of $X$ as a solution to Eq.~\eqref{EqRelControlSDDE} under the relaxed control $(R,W)$ with initial condition $\phi$. Notice that $X$ is continuous on $[0,\infty)$ because of the condition on the jump heights of $(\xi^{M})$, cf.\ Theorem~3.10.2 in \citet[p.\,148]{ethierkurtz86}. Let us define c\`adl\`ag processes $C^{M}$, $C$ on $[0,\infty)$ by
\begin{align*}
    C^{M}(t)&:=\; \phi^{M}(0) \;+\; \int_{0}^{t} b\bigl(\bar{\xi}^{M}_{s},\bar{u}^{M}(s)\bigr)\,ds
\;+\; \varepsilon^{M}_{2}(t), && t \geq 0, & \\
C(t)&:=\; \phi(0) \;+\; \int_{\Gamma\times[0,t]} b(X_{s},\gamma)\, dR(s,\gamma), && t \geq
0. &
\end{align*}
Then $C$, $C^{M}$ are bounded on compact time intervals uniformly in $M \in \mathbb{N}$. Invoking Skorohod's representation theorem, one establishes weak convergence of $(C^{M})$ to $C$ as in the proof of Proposition~\ref{PropRelaxedCompactness}.

The sequence $(\bar{W}^{M})$ is of \emph{uniformly controlled variations}, hence a \emph{good} sequence of integrators in the sense of \citet{kurtzprotter91}, because the jump heights are uniformly bounded and $\bar{W}^{M}$ is a martingale for each $M \in \mathbb{N}$. We have weak convergence of $(\bar{W}^{M})$ to $W$. The results in \citet{kurtzprotter91} guarantee weak convergence of the corresponding adapted quadratic variation processes, that is $([\bar{W}^{M},\bar{W}^{M}])$ converges weakly to $[W,W]$ in $\tilde{D}_{\infty} = D_{\mathbb{R}}([0,\infty))$, where the square brackets indicate the adapted quadratic (co-)variation. Convergence also holds for the sequence of process pairs $(\bar{W}^{M},[\bar{W}^{M},\bar{W}^{M}])$ in $D_{\mathbb{R}^{2}}([0,\infty))$, see Theorem~36 in \citet{kurtzprotter04}.

We now know that each of the sequences $(\bar{\xi}^{M})$, $(C^{M})$, $(\bar{W}^{M})$, $([\bar{W}^{M},\bar{W}^{M}])$ is weakly convergent in $D_{\mathbb{R}}([0,\infty))$. Actually, we have weak convergence for the sequence of process quadruples $(\bar{\xi}^{M},C^{M},\bar{W}^{M},[\bar{W}^{M},\bar{W}^{M}])$ in $D_{\mathbb{R}^{4}}([0,\infty))$. To see this notice that each of the sequences $(\bar{\xi}^{M}+C^{M})$, $(\bar{\xi}^{M}+\bar{W}^{M})$, $(\bar{\xi}^{M}+[\bar{W}^{M},\bar{W}^{M}])$, $(C^{M}+\bar{W}^{M})$, $(C^{M}+[\bar{W}^{M},\bar{W}^{M}])$, and $(\bar{W}^{M}+[\bar{W}^{M},\bar{W}^{M}])$ is tight in $D_{\mathbb{R}}([0,\infty))$, because the limit processes $C$, $X$, $W$, and $[W,W] = Id_{[0,\infty)}$ are all continuous on $[0,\infty)$. According to Problem~22 in \citet[p.\,153]{ethierkurtz86} this implies tightness of the quadruple sequence in $D_{\mathbb{R}^{4}}([0,\infty))$. Since the four component sequences are all weakly convergent, the four-dimensional sequence must have a unique limit point, namely $(X,C,W,[W,W])$. By virtue of Skorohod's theorem, we may again work under $\Prb$-almost sure convergence. Since $C$, $X$, $W$, $[W,W]$ are all continuous, it follows that $C^{M} \to C$, $\bar{\xi}^{M} \to X$, $\bar{W}^{M} \to W$, $[\bar{W}^{M},\bar{W}^{M}] \to [W,W]$ uniformly on compact subintervals of $[0,\infty)$ with probability one.

Define the mapping $F\!: D_{0}\times\tilde{D}_{\infty} \rightarrow \tilde{D}_{\infty}$ by
\begin{equation*}
    F(\phi,x)(t):=\; \sigma\left([-r,0] \ni s \mapsto \begin{cases} x(t\!+\!s) &\text{if}\; t\!+\!s \geq 0, \\ \phi(t\!+\!s) &\text{else}\end{cases} \right),\quad t \geq 0.
\end{equation*}
For $M \in \mathbb{N}$, let $F^{M}$ be the mapping from $\tilde{D}_{\infty}$ to $\tilde{D}_{\infty}$ given by $F^{M}(x)\!:=F(\phi^{M},x)$. Let $H^{M}\!: \tilde{D}_{\infty} \rightarrow \tilde{D}_{\infty}$ be the c\`adl\`ag interpolation operator of degree $M$, that is $H^{M}(x)$ is the piecewise constant c\`adl\`ag interpolation to $x \in \tilde{D}_{\infty}$ along the time grid of mesh size $\frac{r}{M}$ starting at zero. Define $\bar{F}^{M}\!: \tilde{D}_{\infty} \rightarrow \tilde{D}_{\infty}$ by
\begin{equation*}
    \bar{F}^{M}(x)(t):=\; F\bigl(\phi^{M},H^{M}(x)\bigr) \bigl(\lfloor t\rfloor_{M}\bigr), \quad t \geq 0,
\end{equation*}
where $\lfloor t\rfloor_{M}\!:= \frac{r}{M}\lfloor\tfrac{M}{r}t\rfloor$. If $\psi \in D_{\infty}$, we will take $F^{M}(\psi)$, $\bar{F}^{M}(\psi)$ and $F(\psi)$ to equal $F^{M}(x)$, $\bar{F}^{M}(x)$ and $F(\phi,x)$, respectively, where $x$ is the restriction of $\psi$ to $[0,\infty)$. Eq.~\eqref{EqApproximateSDDE} translates to
\begin{equation*}
    \bar{\xi}^{M}(t) \;=\; C^{M}(t) \;+\; \int_{0}^{t} \bar{F}^{M}(\bar{\xi}^{M})(s-)d\bar{W}^{M}(s), \quad t \geq 0.
\end{equation*}
Let $\hat{\xi}$ be the unique c\`adl\`ag process solving
\begin{align*}
    &\hat{\xi}(s) \;=\; \phi(s), \quad s \in [-r,0), & &\hat{\xi}(t) \;=\; C(t) \;+\; \int_{0}^{t} F(\hat{\xi})(s-)dW(s), \quad t \geq 0.&
\end{align*}
Fix $T > 0$. Since $\bar{\xi}^{M}$ converges to $X$ as $M$ goes to infinity uniformly on compacts with probability one, it is enough to show that
\begin{equation} \label{EqPrfSquareError} \tag{$\ast$}
    \Mean\Bigl(\sup_{t\in[-r,T]} \bigl|\hat{\xi}(t)-\bar{\xi}^{M}(t)\bigr|^{2} \Bigr) \quad\stackrel{M\to\infty}{\longrightarrow}\quad 0.
\end{equation}
First observe that
\begin{align*}
    & \Mean\Bigl(\sup_{t\in[0,T]} \bigl|C(t)-C^{M}(t)\bigr|^{2} \Bigr) \quad\stackrel{M\to\infty}{\longrightarrow}\; 0, & & \sup_{t\in[-r,0)} \bigl|\hat{\xi}(t)-\bar{\xi}^{M}(t)\bigr|^{2} \quad\stackrel{M\to\infty}{\longrightarrow}\; 0, &
\end{align*}
because $C$ is uniformly bounded on compact time intervals and
$\phi$ is c\`adl\`ag and continuous on $[-r,0)$. Given $\varepsilon
> 0$, by Lemma~\ref{LemmaIntegralEstimate} in the Appendix and by Gronwall's lemma we find that there is a positive
number $M_{0} = M_{0}(\varepsilon)$ such that for all $M \geq M_{0}$
\begin{equation*}
    \Mean\Bigl(\sup_{t\in[0,T]} \Bigl| \int_{0}^{t}F(\hat{\xi})(s-)dW(s)
- \int_{0}^{t} \bar{F}^{M}(\bar{\xi}^{M})(s-)d\bar{W}^{M}(s)
\Bigr|^{2}\Bigr) \;\leq\; 76T\varepsilon(K^2\!+1)\exp\bigl(4K^{2}_{L}T\bigr).
\end{equation*}
This yields \eqref{EqPrfSquareError} and the assertion follows.
\end{proof}

If we consider approximations along all equidistant partitions of
$[-r,0]$, then the hypothesis about the uniform convergence of the
initial conditions implies that $\phi$ must be continuous on
$[-r,0]\setminus \{0\}$. In case $\phi$ has jumps at
positions locatable on one of the equidistant partitions, the convergence results continue
to hold when we restrict to a sequence of refining partitions.

\section{Convergence of the minimal costs} \label{SectCosts}

The objective behind the introduction of sequences of approximating chains was to obtain a device for approximating the value function $V$ of the original problem. The idea now is to define, for each discretisation degree $M \in \mathbb{N}$, a discrete control problem with cost functional $J^{M}$ so that $J^{M}$ is an approximation of the cost functional $J$ of the original problem in the following sense: Given a suitable initial segment $\phi \in D_{0}$ and a sequence of discrete admissible controls $(u^{M})$ such that $(\bar{u}^{M})$ weakly converges to a relaxed control $R$, we have $J(\phi,u^{M}) \to \hat{J}(\phi,R)$ as $M$ tends to infinity. Under the assumptions introduced above, it will follow that also the value functions associated with the discrete cost functionals converge to the value function of the original problem.

Fix $M \in \mathbb{N}$, and let $h\!:= \frac{r}{M}$. Denote by $\mathcal{U}^{M}_{ad}$ the set of discrete admissible controls of degree $M$. Define the \emph{cost functional of degree $M$} by
\begin{equation} \label{ExCostFunctionalM}
    J^{M}\bigl(\phi,u\bigr) :=\; \Mean\left(\sum_{n=0}^{N_{h}-1} \exp(-\beta nh)\cdot k\bigl(\xi(n),u(n)\bigr)\cdot h \;+\; g\bigl(\xi(N_{h})\bigr)\right),
\end{equation}
where $\phi \in D_{0}$, $u \in \mathcal{U}^{M}_{ad}$ is defined on the stochastic basis $(\Omega,\mathcal{F},(\mathcal{F}_{t}),\Prb)$ and $(\xi(n))$ is a discrete chain of degree $M$ defined according to $p^{M}$ and $u$ with initial condition $\phi$. The discrete \emph{exit time step} $N_{h}$ is given by
\begin{equation} \label{ExEndStepM}
    N_{h} :=\; \min\{n \in \mathbb{N}_{0} \;|\; \xi(n) \notin I_{h}\} \;\wedge\; \lfloor\tfrac{\bar{T}}{h}\rfloor.
\end{equation}
Denote by $\bar{\tau}^{M}\!:= h\cdot N_{h}$ the exit time for the corresponding interpolated processes. The \emph{value function of degree $M$} is defined as
\begin{equation} \label{ExValueFunctionM}
    V^{M}(\phi) :=\; \inf\bigl\{J^{M}\bigl(\phi,u\bigr) \;\big|\; u \in \mathcal{U}^{M}_{ad}\bigr\}, \quad \phi \in D_{0}.
\end{equation}
We are now in a position to state the result about convergence of
the minimal costs. Proposition~\ref{PropMinimalCosts} and
Theorem~\ref{ThValueConvergence} are comparable to Theorems 10.5.1
and 10.5.2 in \citet[pp.\,292-295]{kushnerdupuis01}. Let us suppose
that the initial condition $\phi \in D_{0}$ and the sequence of
partitions of $[-r,0]$ are such that the discretised initial
conditions converge to $\phi$ uniformly on $[-r,0]$.

\begin{prop} \label{PropMinimalCosts}
Assume \hypref{HypCadlag}\,--\,\hypref{HypEllipticity}. If the
sequence $(\bar{\xi}^{M},\bar{u}^{M},\bar{W}^{M},\bar{\tau}^{M})$ of
interpolated processes converges weakly to a limit point
$(X,R,W,\tau)$, then $X$ is a solution to
Eq.~\eqref{EqRelControlSDDE} under relaxed control $(R,W)$ with
initial condition $\phi$, $\tau$ is the exit time for $X$ as given
by \eqref{ExEndTime}, and we have
    \[ J^{M}(\phi,u^{M}) \;\stackrel{M\to\infty}{\longrightarrow}\; \hat{J}(\phi,R). \]
\end{prop}

\begin{proof}
The convergence assertion for the costs is a consequence of Proposition~\ref{PropChainConvergence}, the fact that, by virtue of Assumption~\hypref{HypEllipticity}, the exit time $\hat{\tau}$ defined in \eqref{ExExitTime} is Skorohod-continuous, and the definition of $J^{M}$ and $J$ (or $\hat{J}$).
\end{proof}

\begin{thrm} \label{ThValueConvergence}
Assume \hypref{HypCadlag}\,--\,\hypref{HypEllipticity}. Then we have
$\lim_{M\to\infty}V^{M}(\phi)=V(\phi)$.
\end{thrm}

\begin{proof}
First notice that $\liminf_{M\to\infty} V^{M}(\phi) \geq V(\phi)$ as a consequence of Propositions \ref{PropChainConvergence} and \ref{PropMinimalCosts}. In order to show $\limsup_{M\to\infty} V^{M}(\phi) \leq V(\phi)$ choose a relaxed control $(R,W)$ so that $\hat{J}(\phi,R) = V(\phi)$ according to Proposition~\ref{PropRelaxedCompactness}. Given $\varepsilon > 0$, one can construct a sequence of discrete admissible controls $(u^{M})$ such that $((\bar{\xi}^{M},\bar{u}^{M},\bar{W}^{M},\bar{\tau}^{M}))$ is weakly convergent, where $(\bar{\xi}^{M})$, $(\bar{W}^{M})$, $(\bar{\tau}^{M})$ are constructed as above, and
    \[ \limsup_{M\to\infty} |J^{M}(\phi,u^{M})-\hat{J}(\phi,R)| \leq \varepsilon.\]
The existence of such a sequence of discrete admissible controls is guaranteed, cf.\ the discussion at the end of Section~\ref{SectExistence}. By definition, $V^{M}(\phi) \leq J^{M}(\phi,u^{M})$ for each $M \in \mathbb{N}$. Using Proposition~\ref{PropMinimalCosts} we find that
\begin{equation*}
    \limsup_{M\to\infty} V^{M}(\phi) \;\leq\; \limsup_{M\to\infty} J^{M}(\phi,u^{M}) \;\leq\; V(\phi) + \varepsilon,
\end{equation*}
and since $\varepsilon$ was arbitrary, the assertion follows.
\end{proof}

\begin{appendix}
\section*{Appendix}
The proof of the following lemma makes use of standard techniques. In the context of approximation of SDDEs, it should be compared to Section~7 in \citet{mao03}.

\begin{lemma} \label{LemmaIntegralEstimate}
In the notation and under the assumptions of Proposition~\ref{PropChainConvergence} it holds that for every $\varepsilon > 0$ there is $M_{0} \in \mathbb{N}$ such that for all $M \geq M_{0}$
\begin{equation*}
    \begin{split}
    & \Mean\Bigl(\sup_{t\in[0,T]} \bigl| \int_{0}^{t}F(\hat{\xi})(s-)dW(s) - \int_{0}^{t} \bar{F}^{M}(\bar{\xi}^{M})(s-)d\bar{W}^{M}(s) \bigr|^{2}\Bigr) \\[1ex]
    \leq\quad   & 4K^{2}_{L}\, \int_{0}^{T} \Mean\Bigl(\sup_{t\in[-r,s]} \bigl|\hat{\xi}(t)-\bar{\xi}^{M}(t)\bigr|^{2}\Bigr)\,ds \;+\; 76T\varepsilon(K^2+1).
    \end{split}
\end{equation*}
\end{lemma}

\begin{proof} Clearly,
\begin{equation} \label{EqPrfFirstSplitting}
    \begin{split}
    & \Mean\Bigl(\sup_{t\in[0,T]} \bigl| \int_{0}^{t}F(\hat{\xi})(s-)dW(s) - \int_{0}^{t} \bar{F}^{M}(\bar{\xi}^{M})(s-)d\bar{W}^{M}(s) \bigr|^{2}\Bigr) \\[1ex]
    \leq\quad & 2\Mean\Bigl(\sup_{t\in[0,T]} \bigl| \int_{0}^{t}F(\hat{\xi})(s-)dW(s) - \int_{0}^{t}\bar{F}^{M}(\bar{\xi}^{M})(s-)dW(s) \bigr|^{2}\Bigr) \\
    +\; & 2\Mean\Bigl(\sup_{t\in[0,T]} \bigl| \int_{0}^{t}\bar{F}^{M}(\bar{\xi}^{M})(s-)dW(s) - \int_{0}^{t}\bar{F}^{M}(\bar{\xi}^{M})(s-)d\bar{W}^{M}(s) \bigr|^{2}\Bigr)
    \end{split}
\end{equation}
Using Doob's maximal inequality, It{\^o}'s isometry, Fubini's theorem and Assumption~\hypref{HypLipschitz}, for the first expectation on the right hand side of \eqref{EqPrfFirstSplitting} we obtain the estimate
\begin{equation} \label{EqPrfFirstEstimate}
    \begin{split}
    & \Mean\Bigl(\sup_{t\in[0,T]} \bigl| \int_{0}^{t}F(\hat{\xi})(s-)dW(s) - \int_{0}^{t}\bar{F}^{M}(\bar{\xi}^{M})(s-)dW(s) \bigr|^{2}\Bigr) \\[1ex]
    \leq\quad & 4\int_{0}^{T} \Mean\Bigl(\bigl|F(\hat{\xi})(s) - \bar{F}^{M}(\bar{\xi}^{M})(s)\bigr|^{2}\Bigr)\,ds \\[1ex]
    \leq\quad & 4K^{2}_{L}\, \int_{0}^{T} \Mean\Bigl(\sup_{t\in[-r,s]} \bigl|\hat{\xi}(t)-\bar{\xi}^{M}(t)\bigr|^{2}\Bigr)\,ds.
    \end{split}
\end{equation}
Fix any $N \in \mathbb{N}$. The second expectation on the right hand side of \eqref{EqPrfFirstSplitting} splits up into three terms according to
\begin{equation} \label{EqPrfSecondSplitting}
    \begin{split}
    & \Mean\Bigl(\sup_{t\in[0,T]} \bigl| \int_{0}^{t}\bar{F}^{M}(\bar{\xi}^{M})(s-)dW(s) - \int_{0}^{t}\bar{F}^{M}(\bar{\xi}^{M})(s-)d\bar{W}^{M}(s) \bigr|^{2}\Bigr) \\[1ex]
    \leq\quad & 4\Mean\Bigl(\sup_{t\in[0,T]} \bigl| \int_{0}^{t}\bar{F}^{M}(\bar{\xi}^{M})(s-)dW(s) - \int_{0}^{t}\bar{F}^{N}(\bar{\xi}^{M})(s-)dW(s) \bigr|^{2}\Bigr) \\
    +\; & 4\Mean\Bigl(\sup_{t\in[0,T]} \bigl| \int_{0}^{t}\bar{F}^{N}(\bar{\xi}^{M})(s-)dW(s) - \int_{0}^{t}\bar{F}^{N}(\bar{\xi}^{M})(s-)d\bar{W}^{M}(s) \bigr|^{2}\Bigr) \\
    +\; & 4\Mean\Bigl(\sup_{t\in[0,T]} \bigl| \int_{0}^{t}\bar{F}^{N}(\bar{\xi}^{M})(s-)d\bar{W}^{M}(s) - \int_{0}^{t}\bar{F}^{M}(\bar{\xi}^{M})(s-)d\bar{W}^{M}(s) \bigr|^{2}\Bigr).
    \end{split}
\end{equation}
Again using Doob's maximal inequality and a generalized version of It{\^o}'s isometry \citep[cf.][pp.\,73-77]{protter03}, for the first and third expectation on the right hand side of \eqref{EqPrfSecondSplitting} we get
\begin{align}
    \label{EqPrfSecondEstimate} \begin{split}
    & \Mean\Bigl(\sup_{t\in[0,T]} \bigl| \int_{0}^{t}\bar{F}^{M}(\bar{\xi}^{M})(s-)dW(s) - \int_{0}^{t}\bar{F}^{N}(\bar{\xi}^{M})(s-)dW(s) \bigr|^{2}\Bigr) \\[1ex]
    \leq\quad & 4\Mean\Bigl(\int_{0}^{T} \bigl|\bar{F}^{M}(\bar{\xi}^{M})(s) - \bar{F}^{N}(\bar{\xi}^{M})(s)\bigr|^{2}\,ds\Bigr)
    \end{split} \\
    \intertext{and}
    \label{EqPrfThirdEstimate} \begin{split}
    & \Mean\Bigl(\sup_{t\in[0,T]} \bigl| \int_{0}^{t}\bar{F}^{N}(\bar{\xi}^{M})(s-)d\bar{W}^{M}(s) - \int_{0}^{t}\bar{F}^{M}(\bar{\xi}^{M})(s-)d\bar{W}^{M}(s) \bigr|^{2}\Bigr) \\[1ex]
    \leq\quad & 4\Mean\Bigl(\int_{0}^{T} \bigl|\bar{F}^{M}(\bar{\xi}^{M})(s-) - \bar{F}^{N}(\bar{\xi}^{M})(s-)\bigr|^{2}\,d\bigl[\bar{W}^{M},\bar{W}^{M}\bigr](s)\Bigr)
    \end{split}
\end{align}
respectively. Notice that, path-by-path, we have
\begin{equation*}
    \begin{split}
     & \int_{0}^{T} \bigl|\bar{F}^{M}(\bar{\xi}^{M})(s-) - \bar{F}^{N}(\bar{\xi}^{M})(s-)\bigr|^{2}\,d\bigl[\bar{W}^{M},\bar{W}^{M}\bigr](s) \\[1ex]
    \leq\quad &  \sum_{i=0}^{\lfloor\tfrac{M}{r}T\rfloor} \bigl|\bar{F}^{M}(\bar{\xi}^{M})\bigl(\tfrac{r}{M}i\bigr) - \bar{F}^{N}(\bar{\xi}^{M})\bigl(\tfrac{r}{M}i\bigr)\bigr|^{2}\cdot \Bigl([\bar{W}^{M},\bar{W}^{M}\bigr]\bigl(\tfrac{r}{M}(i\!+\!1)\bigr)-[\bar{W}^{M},\bar{W}^{M}\bigr]\bigl(\tfrac{r}{M}i\bigr)\Bigr).
    \end{split}
\end{equation*}
In order to estimate the second expectation on the right hand side of \eqref{EqPrfSecondSplitting}, observe that, $\Prb$-almost surely, for all $t \in [0,T]$
\begin{equation*}
    \begin{split}
    \int_{0}^{t}\bar{F}^{N}(\bar{\xi}^{M})(s-)dW(s) \quad=\quad & \bar{F}^{N}(\bar{\xi}^{M})\bigl(\lfloor t\rfloor_{N}\bigr)\cdot \Bigl(W(t)-W\bigl(\lfloor t\rfloor_{N}\bigr)\Bigr) \\
    +\; & \sum_{i=0}^{\lfloor\tfrac{N}{r}t\rfloor-1} \bar{F}^{N}(\bar{\xi}^{M})\bigl(\tfrac{r}{N}i\bigr)\cdot \Bigl(W\bigl(\tfrac{r}{N}(i\!+\!1)\bigr)-W\bigl(\tfrac{r}{N}i\bigr)\Bigr),
    \end{split}
\end{equation*}
as $F^{N}(\bar{\xi}^{M})$ is piecewise constant on the grid of mesh size $\frac{r}{N}$. On the other hand,
\begin{equation*}
    \begin{split}
    \int_{0}^{t}\bar{F}^{N}(\bar{\xi}^{M})(s-)d\bar{W}^{M}(s) \quad=\quad & \bar{F}^{N}(\bar{\xi}^{M})\bigl(\lfloor t\rfloor_{N}\bigr)\cdot \Bigl(\bar{W}^{M}(t)-\bar{W}^{M}\bigl(\lfloor t\rfloor_{N}\bigr)\Bigr) \\
    +\; & \sum_{i=0}^{\lfloor\tfrac{N}{r}t\rfloor-1} \bar{F}^{N}(\bar{\xi}^{M})\bigl(\tfrac{r}{N}i\bigr)\cdot \Bigl(\bar{W}^{M}\bigl(\tfrac{r}{N}(i\!+\!1)\bigr)-\bar{W}^{M}\bigl(\tfrac{r}{N}i\bigr)\Bigr).
    \end{split}
\end{equation*}
By Assumption~\hypref{HypGrowth}, $|\sigma|$ is bounded by a constant $K$, hence
\begin{equation*}
    \begin{split}
    & \bigl| \int_{0}^{t}\bar{F}^{N}(\bar{\xi}^{M})(s-)dW(s) - \int_{0}^{t}\bar{F}^{N}(\bar{\xi}^{M})(s-)d\bar{W}^{M}(s) \bigr| \\[1ex]
    \leq\quad & 2K\lfloor\tfrac{N}{r}t\rfloor\cdot\! \sup_{s\in[0,t]} |W(s)-\bar{W}^{M}(s)| \quad\leq\quad 2K\tfrac{N}{r}T\cdot\! \sup_{s\in[0,T]} |W(s)-\bar{W}^{M}(s)|.
    \end{split}
\end{equation*}
Bounded convergence yields for each fixed $N \in \mathbb{N}$
\begin{equation} \label{EqPrfMiscIntConvergence}
    \Mean\Bigl(\sup_{t\in[0,T]} \bigl| \int_{0}^{t}\bar{F}^{N}(\bar{\xi}^{M})(s-)dW(s) - \int_{0}^{t}\bar{F}^{N}(\bar{\xi}^{M})(s-)d\bar{W}^{M}(s) \bigr|^{2}\Bigr) \quad\stackrel{M\to\infty}{\longrightarrow} 0.
\end{equation}
Let $x, y \in \tilde{D}_{\infty}$. By Assumption~\hypref{HypLipschitz} we have for all $t \in [0,T]$
\begin{equation*}
    \begin{split}
    & \bigl|\bar{F}^{N}(y)(t) - F(\phi,x)(t)\bigr| \quad=\quad \bigl|F\bigl(\phi^{N},H^{N}(y)\bigr)(\lfloor t\rfloor_{N}) - F(\phi,x)(t)\bigr| \\[1ex]
    \leq\quad & K_{L}\cdot\! \sup_{s\in[-r,0)}\bigl|\phi^{N}(s)-\phi(s)\bigr| \;+\; K_{L}\cdot\! \sup_{s\in[0,T]}\bigl|H^{N}(y)(s)-x(s)\bigr| \\
    +\; & \bigl|F(\phi,x)\bigl(\lfloor t\rfloor_{N}\bigr) - F(\phi,x)(t)\bigr|.
    \end{split}
\end{equation*}
By Assumption~\hypref{HypCadlag}, the map $[0,T] \ni t \mapsto F(\phi,x)(t)$ is c\`adl\`ag, whence it has only finitely many jumps larger than any given positive lower bound. Thus, given $\varepsilon > 0$, there is a finite subset $A = A(\varepsilon,T,\phi,x) \subset [0,T]$ such that
\begin{equation*}
    \limsup_{N\to\infty}\; \bigl|F(\phi,x)\bigl(\lfloor t\rfloor_{N}\bigr) - F(\phi,x)(t)\bigr| \quad\leq\; \varepsilon \quad\text{for all}\; t \in [0,T] \setminus A.
\end{equation*}
Moreover, the convergence is uniform in the following sense \citep[cf.][]{billingsley99}: We can choose the finite set $A$ in such a way that there is $N_{0} = N_{0}(\varepsilon,T,\phi,x) \in \mathbb{N}$ so that
\begin{equation*}
    \bigl|F(\phi,x)\bigl(\lfloor t\rfloor_{N}\bigr) - F(\phi,x)(t)\bigr| \quad\leq\; 2\varepsilon \quad\text{for all}\; t \in [0,T] \setminus A,\; N \geq N_{0}.
\end{equation*}
Given $\varepsilon > 0$, we therefore find $N \in \mathbb{N}$ and an event $\tilde{\Omega}$ with $\Prb(\tilde{\Omega}) \geq 1\!-\!\varepsilon$ so that for each $\omega \in \tilde{\Omega}$ there is a finite subset $A_{\omega} \subset [0,T]$ with $\#A_{\omega} \leq N\varepsilon$ and such that for all $t \in [0,T] \setminus A_{\omega}$ and all $M \geq N$ we have
\begin{equation*}
    \bigl|\bar{F}^{M}\bigl(\bar{\xi}^{M}(\omega)\bigr)(t) - F\bigl(X(\omega)\bigr)(t)\bigr|^{2} \;+\; \bigl|\bar{F}^{N}\bigl(\bar{\xi}^{M}(\omega)\bigr)(t) - F\bigl(X(\omega)\bigr)(t)\bigr|^{2} \quad\leq\quad \varepsilon.
\end{equation*}
The expression on the right hand side of \eqref{EqPrfSecondEstimate} is then bounded from above by $9T\varepsilon(K^{2}+1)$. For $M$ big enough, also the expression on the right hand side of \eqref{EqPrfThirdEstimate} is smaller than $9T\varepsilon(K^{2}+1)$, and the expectation in \eqref{EqPrfMiscIntConvergence} is smaller than $T\varepsilon$.
\end{proof}

\end{appendix}

\newpage


\bibliographystyle{plainnat}
\bibliography{MathLit}

\end{document}